\documentclass{elsart3-1-bis}
\usepackage{amssymb}
\usepackage[english,francais]{babel}

\newtheorem{theorem}{Theorem}[section]

\newtheorem{e-proposition}[theorem]{Proposition}
\newtheorem{corollary}[theorem]{Corollary}
\newtheorem{e-definition}[theorem]{Definition\rm}

\newtheorem{theoreme}{Th{\'e}or{\`e}me}[section]

\renewcommand{\theequation}{\arabic{equation}}
\def\og{\leavevmode\raise.3ex\hbox{$\scriptscriptstyle\langle\!\langle$~}}
\def\fg{\leavevmode\raise.3ex\hbox{~$\!\scriptscriptstyle\,\rangle\!\rangle$}}

\def \a {\alpha}
\def \p {{\mathbb P}}
\def \C {{\mathbb C}}

\def \G {{\mathbb G}}

\def \Z {{\mathbb{Z}}}

\def \pic {{\rm Pic}}
\def \codim {{\rm codim}}
\def \noi {\noindent}
\def \vs {\vskip}
\def \hs {\hskip}
\def \hod {{\rm Hodge}}
\newtheorem{question}[theoreme]{\it Question\/}
\newenvironment{preu}{\noi\sc{Proof}\hs 0.1 cm 
--- \rm }{\hfill$\Box$\vs 0.2 cm}

\begin{document}

~

\vs -4 cm
\begin{frontmatter}

\selectlanguage{english}
\title{Small codimension smooth subvarieties in even-dimensional
  homogeneous spaces with Picard group $\Z$}

\selectlanguage{english}
\author[authorlabel1]{Nicolas Perrin}
\ead{nperrin@math.jussieu.fr}

\address[authorlabel1]{Universit{\'e} Pierre et Marie Curie - Paris 6\\
Institut de Math{\'e}matiques de Jussieu\\
175 rue du Chevaleret, 75013 Paris.}

\vs 1 cm


\begin{abstract}
\selectlanguage{english}
We investigate a method proposed by E. Arrondo and J. Caravantes to study
  the Picard group of a smooth low-codimension subvariety $X$ in a
  variety $Y$ when $Y$ is homogeneous. We prove that this method is
  strongly related to the signature $\sigma_Y$ of the Poincar{\'e} pairing
  on the middle cohomology of $Y$. We give under some topological
  assumptions a bound on the rank of Picard group $\pic(X)$ in terms of
  $\sigma_Y$ and remove these assumptions for grassmannians to recover
  the main result of E. Arrondo and J. Caravantes.
\end{abstract}
\end{frontmatter}

\selectlanguage{english}

\section{Introduction}

E. Arrondo and J. Caravantes present in \cite{arrondo} a
sketch of program to decide if a divisor $D$ in a smooth subvariety
$X$ of a smooth even-dimensional variety $Y$ is 
equivalent to a multiple of an hyperplane section of $X$. 
To apply this method, they need some hard
topological results that we will not discuss here.
We will however discuss the general validity of the rest of the method
in the case of homogeneous varieties $Y$ with Picard group $\Z$ and
prove that
it is closely related to the signature $\sigma_Y$ of the Poincar{\'e}
pairing on the middle cohomology of $Y$.
Let us fix the following notations: $N=\dim(Y)$ the dimension of $Y$
is even and to keep coherent notations with \cite{arrondo} we set
$N=2(n-1)$. Recall the definition (see \cite{De}):

\vs 0.2 cm

\begin{e-definition}
  A subvariety $X$ of $Y$ is said to be \emph{cumbersome} if for
  any subvariety $Z\subset Y$ such that $\dim(Z)\geq \codim(X)$ we
  have $[X]\cdot[Z]\neq0$.
\end{e-definition}

\vs 0.2 cm

We will prove the following:

\vs 0.2 cm

\begin{theorem}
  \label{main}
Let 
$X$ be a cumbersome
smooth subvariety in $Y$ of dimension $n'$ with $n'\geq n$.
Assume that $X$ is simply connected and that, for any $n$-dimensional
Schubert subvariety $Y'$ of $Y$, the intersection of $X$ with a general
translate of $Y'$ is irreducible, then $\pic(X)$ 
is free
of rank at most
$
{\frac{1}{2}(h^{n-1}(Y,\C)-\sigma_Y)+1.}$ 
\end{theorem}

\vs 0.2 cm

One major difficulty to use this result is to prove the two
topological assumptions on $X$, namely that $X$ is simply connected
and that its intersection with the general translate of a Schubert
subvariety of dimension $n$ is irreducible. However,
in the case where $Y$ is the Grassmannian variety $\G(p,m)$,
O. Debarre \cite{De} proved such results. In that case, set
$R=\max(p,m-p)$, $r=\min(p,m-p)$, $\delta=R-1$ for $r\geq3$ and
$\delta=1$ for $r\leq2$. We deduce the following:

\vs 0.2 cm

\begin{corollary}
  \label{gras}
Assume $X$ is a smooth cumbersome subvariety of dimension
$n'\geq n-1+\delta$ of $Y$ then $\pic(X)$ is free of rank at most
$
{\frac{1}{2}(h^{n-1}(Y,\C)-\sigma_Y)+1.}$
\end{corollary}

\vs 0.2 cm

\noi
Finally we recover the main result in \cite{arrondo}:

\vs 0.2 cm

\begin{corollary}
  \label{gras-lin}
Assume that $r\leq2$,
or equivalently that $Y$ is a projective
space or a Grassmannians of lines, 
then any smooth cumbersome
subvariety $X$ of $Y$ of dimension $n'\geq n$ satisfies $\pic(X)=\Z$.
\end{corollary}

\vs 0.2 cm

Remark that stronger results were proved by A. J. Sommese \cite{So}
but in smaller codimension.
In this note we prove that the method of E. Arrondo and J. Caravantes
to show that the Picard group of $X$ is $\Z$ will not work in
general but only if the Poincar{\'e} duality is positive definite on the
middle cohomology. The only\footnote{See \cite{Sl}}
even-dimensional homogeneous varieties with Picard group $\Z$ and
positive definite Poincar{\'e} pairing on the middle 
cohomology are the projective spaces, the grassmannians of lines, the
quadrics of dimension multiple by 4 and the Cayley plane
$\mathbb{O}\p^2$. It is striking that these are examples of projective
spaces over composition algebras (over $\mathbb{C}$, $\mathbb{H}$ and
$\mathbb{O}$, cf. \cite{Chaput}). With this in mind we ask the
following: 

\vs 0.2 cm

  \begin{question}
    Is it true that, for any smooth cumbersome subvariety $X$ of
    dimension at least 9 of the Cayley plane $\mathbb{O}\p^2$, we have
    $\pic(X)=\Z$~?
  \end{question}

\vs 0.2 cm

\noi
{\bf Acknowledgements:} I would like to thank L. Gruson for many
usefull discussion during the preparation of this text.

\section{Proof of the Theorem}

To prove the theorem we will in some sense restrict ourselves to the
$n$-dimensional case. 
Let $X'$ be a variety satisfying the
hypothesis of Theorem \ref{main}. Let $D'$ be a divisor of $X'$, we
may assume $D'$ to be smooth\footnote{This will be harmless by
  replacing $D'$ by $D'+mH_{X'}$ with $m$ big.}. We will denote by $X$ a
general hyperplane section of $X'$ of dimension $n$ and by $D$ the
divisor in $X$ corresponding to $D'$. Remark that $X$ and $D$ are also
smooth and that $X$ is cumbersome.

\vs 0.2 cm

It is a classical result that a basis of $H^*(Y,\Z)$ is given by the
classes $\sigma(w)$ of Schubert subvarieties where $w$ is in a coset
$W_Y$ of the Weyl group $W$ (see for example \cite{Brion}). We will
denote by $(\sigma(t))_{t\in T}$ and $(\sigma(u))_{u\in U}$ the
classes of Schubert subvarieties of dimensions $n$ and $n-1$. They
form a basis of the effective monoids in $H^{n-2}(Y,\Z)$ and
$H^{n-1}(Y,\Z)$. We use the notation:
${\sigma(w)\cdot\sigma(w')=
  \sum_{w''}c_{w,w'}^{w''}\sigma(w'')}$ and write
$$\displaystyle{[X]=\sum_{t\in T}a_t\sigma(t)\ \ \ {\rm and}\ \ \ 
  [D]=\sum_{u\in U}\a_u\sigma(u)}$$
with $a_t$ and $\a_u$ non negative for all $t\in T$ and $u\in U$
(because $X$ is cumbersome, all the $a_t$ are positive).

\vs 0.2 cm

Poincar{\'e} duality acts as a permutation between Schubert
classes and induces an involution 
$w\mapsto w^*$ on $W_Y$.
In particular, it stabilises $U$ and we may consider the partition
$U_{id}\coprod U_{hyp}=U$ where $U_{id}$ is the set of fixed points of
$U$ under the Poincar{\'e} involution. This 
induces a decomposition of the abelian group $H^{n-1}(Y,\Z)$ into
$H^{n-1}_{id}(Y,\Z)\oplus H^{n-1}_{hyp}(Y,\Z)$. Poincar{\'e} duality acts
as the identity on $H^{n-1}_{id}(Y,\Z)$ and is hyperbolic on
$H^{n-1}_{hyp}(Y,\Z)$. In particular if $\sigma_Y$ is its signature on
$H^{n-1}(Y,\Z)$, we have the equalities 
$\dim(H^{n-1}_{id}(Y,\Z))=\sigma_Y\ \ {\rm
  and}\ \ \dim(H^{n-1}_{hyp}(Y,\Z))=h^{n-1}(Y,\Z)-\sigma_Y.$

\vs 0.2 cm

We first 
compare $D$ and $H_X$ with respect to numerical
equivalence. Consider the 
classes of curves
$C(u)=\sigma(u)\cdot[X]\ {\rm and}\ C(t)=[D]\cdot \sigma(t)$.
We define some
non negative integers by:
$x_u:=H_X\cdot C(u),\ $ $y_u:=D\cdot C(u) ,\ z_t:=H_X\cdot C(t) \ {\rm
  and}\ \lambda_t:=[D]\cdot C(t)$. 
We have the equalities:
$${x_u=\sum_{t\in T}c_{t,H}^{u^*}a_t},\ \ y_u=\a_{u^*}\ \ {\rm and}\ \
{z_t=\sum_{u\in U}c_{t,H}^{u^*}\a_{u}.}$$
Define a matrix $M$ with two lines and
columns indexed by $U\cup T$ by 
$$\displaystyle{M=(m_{i,j})_{i\in\{1,2\},\ j\in U\cup T}=
  \left(
\begin{array}{cc}
(x_u)_{u\in U}& (z_t)_{t\in T}\\
(y_u)_{u\in U}& (\lambda_t)_{t\in T}\\
\end{array}\right)}.$$
The same proof as in \cite{arrondo} and Lefschetz's hyperplane Theorem
lead to the following\footnote{We can apply a numerical version of
  Leschetz Theorem because $X'$ and all its linear sections are simply
  connected thanks to the homotopic Lefschetz's hyperplane Theorem,
  see \cite[Theorem 3.1.21]{La}}: 

\vs 0.2 cm

\begin{e-proposition}
\label{rang1}
  The divisor $D$ is numerically equivalent to a multiple of $H_X$ iff
$M$ is of rank one.
\end{e-proposition}

\vs 0.2 cm

\begin{preu}
We follow the proof of E. Arrondo and J. Caravantes in \cite{arrondo}.
  If $D$ is numerically equivalent to a multiple of $H_X$, then the
  rank of the matrix $M$ is one. Conversely, if the second line of the
  matrix is $q$ times the first one, let $D_0=[D]-qH_X$. The
  intersections $D_0\cdot C(u)$ and $D_0\cdot C(t)$
  vanish. But $H_X^{n-1}$ and $H_X^{n-2}$ are linear combinations of
  $C(u)$ and $C(t)$ respectively, hence if $S$ is
  a smooth surface obtained from $X$ by $n-2$ hyperplane sections, we
  get that 
$$D_0\vert_SH\vert_S=D_0H_X^{n-1}=0\ \ {\rm and}\ \
D_0\vert_S^2=D_0(D-qH_X)H_X^{n-2}=0.$$
By Hogde index Theorem, the divisor $D_0\vert_S$ has to be numerically
trivial. By Lefschetz's hyperplane Theorem, this has to be true for
$D_0$.
\end{preu}

\vs 0.2 cm

\begin{corollary}
\label{rang1-bis}
  The divisor $D'$ is numerically equivalent to a multiple of $H_{X'}$
  iff 
$M$ is of rank one.
\end{corollary}

\vs 0.2 cm

\begin{preu}
  We know from the previous Proposition that $D$ and $H_X$ are
  numerically dependent. Lefschetz's hyperplane theorem implies that
  it is also tha case for $D'$ and $H_{X'}$.
\end{preu}

As in \cite{arrondo}, the sequence $0\to N_{D/X}\to N_{D/Y}\to
(N_{X/Y})\vert_D\to 0$ is exact because $X$ and $D$ are smooth. Taking
the top Chern classes gives the equality
$P:=D\cdot_Y D-(D\cdot_X D)\cdot_X X\vert_X=0$.
In terms of the variables $(a_t)_{t\in T}$, $(\a_u)_{u\in U}$ and
$(\lambda_t)_{t \in T, v\in V}$, we get:
$$P=\sum_{u\in U}\a_u\a_{u^*}-\sum_{t\in T}a_t\lambda_t.$$

The next step is the elimination of the variables $(\lambda_t)_{t\in
  T, v\in V}$. We consider for this, following \cite{arrondo}, the
surfaces $(S_t)_{t\in T}$ where $S_t$ is the intersection of $X$ with
a general subvariety of class $\sigma(t)$. We set:
$$\displaystyle{\hod_t:=\left\vert
  \begin{array}{cc}
H_{S_t}^2&H_{S_t}D_{S_t}\\
H_{S_t}D_{S_t}&D_{S_t}^2
  \end{array}
\right\vert=\left\vert
  \begin{array}{cc}
\displaystyle{\sum_{u\in U}c_{t,H}^ux_u}& z_t\\
\displaystyle{\sum_{u\in U}c_{t,H}^uy_u}
&\lambda_t
  \end{array}
\right\vert}=\left\vert
  \begin{array}{cc}
\displaystyle{\sum_{u\in U}c_{t,H}^ux_u}& \displaystyle{\sum_{u\in
U}c_{t,H}^uy_u}\\ 
\displaystyle{\sum_{u\in U}c_{t,H}^uy_u}
&\lambda_t
  \end{array}
\right\vert=\sum_{u\in U}c_{t,H}^u\left\vert
  \begin{array}{cc}
m_{1,u}& m_{1,t}\\ 
m_{2,u}& m_{2,t}
  \end{array}
\right\vert.$$
Observe that $\hod_t$ is a non negative
linear combination of some minors of the matrix $M$. By hypothesis on
$X'$, the surface $S_t$ is irreducible hence $\hod_t$ is non positive. 
Set $d_t:=H_{S_t}^2>0$
and define the following quadratic form $q$ in the variables
$(\a_u)_{u\in U}$ (observe that $q$ is non positive)
\begin{equation}
  \label{q}
\hs 5.3 cm q:=P+\sum_{t\in T}\frac{a_t}{d_t}\hod_t.
\end{equation}

\vs 0.1 cm

In \cite{arrondo}, E. Arrondo and J. Caravantes decompose $q$, in the
special case of grassmannians of lines, 
as a sum of squares to prove its non negativity. We do this in full
generality and show that some negative part may appear. We know that
all the variables $(x_u)_{u\in U}$ do no vanish because $X$ is
cumbersome. Hence we may define two 
quadratic forms $q'$ and $q''$ in the variables $(\a_u)_{u\in U}$ by:
$$q'=\frac{1}{2}\sum_{t\in T}\sum_{(u,u')\in U^2}
\frac{a_tc_{t,H}^uc_{t,H}^{u'}}{d_tx_{u}x_{u'}} 
\left\vert
\begin{array}{cc}
m_{1,u}&m_{1,u'}\\
m_{2,u}&m_{2,u'}
\end{array}
\right\vert^2\ \ \ \ {\rm and}\ \ \ \ q''=\frac{1}{2}\sum_{u\in
  U}\frac{1}{x_{u}x_{u^*}}\left\vert 
\begin{array}{cc}
  m_{1,u}&m_{1,u^*}\\
m_{2,u}&m_{2,u^*}
\end{array}
\right\vert^2.$$

\vs 0.2 cm

\begin{e-proposition}
\label{formule}
  The formula $q=q'-q''$ holds.
\end{e-proposition}

\vs 0.2 cm

\begin{preu}
For any couple $(u,u')$ of elements in $U$, set 
$\displaystyle{A_{u,u'}=-\sum_{t\in T} \frac{a_t}{d_t}c_{t,H}^{u^*}
  c_{t,H}^{{u'}^*}}.$

The coefficient of $\a_{u}\a_{u'}$ in $q$ is equal to
$A_{u,u'}+A_{u',u}=2A_{u,u'}$ if $u'\neq u^*$ and $u'\neq u$, it is
equal to $2+2A_{u,u'}$ if $u'=u^*\neq u$, to $A_{u,u'}$ if $u^*\neq u'=u$ and to $1+A_{u,u'}$ if $u'=u^*=u$.
The coefficient of $\a_{u}\a_{u'}$ in $q'$ is equal to $2A_{u,u'}$ if $u'\neq u$ and, if $u'=u$, to
$$\displaystyle{\sum_{t\in T}\sum_{u''\in U,\ u''\neq
    u}\frac{a_t}{d_t}c_{t,H}^{u^*}c_{t,H}^{{u''}^*}
  \frac{x_{{u''}^*}}{x_{u^*}} = \sum_{t\in T}\frac{a_t}{d_tx_{u^*}}
  c_{t,H}^{u^*}(d_t-c_{t,H}^{u^*}x_{u^*})} =\frac{x_u}{x_{u^*}}-
\sum_{t\in T}\frac{a_t}{d_t}{c_{t,h}^{u^*}}^2=
\frac{x_u}{x_{u^*}}+A_{u,u'}.$$
The non vanishing coefficients of $\a_{u}\a_{u'}$ in $q''$ are
equal to $-2$ if $u'=u^*\neq u$ and to 
$
{\frac{x_u}{x_{u^*}}} \textrm{
  if }u'=u\neq u^*.$  
But this is exactely the difference $q'-q$.
\end{preu}

\vs 0.2 cm

\begin{corollary}
  \label{multi}
  If for all $u\in U$, we have 
$\displaystyle{\left\vert 
\begin{array}{cc}
  m_{1,u}&m_{1,u^*}\\
m_{2,u}&m_{2,u^*}
\end{array}
\right\vert}=0$,  
then $D'$ is a multiple of $H_{X'}$ in $\pic(X')$.
\end{corollary}

\vs 0.2 cm

\begin{preu}
  Indeed, in that case $q''$ vanishes so $q$ must vanish (and also
  $q'$). Because $X$ is cumbersome, all $a_t$ are different from 0 and
  this implies that all $2\times 2$ minors of $M$ vanish and 
  $M$ is of rank one. Apply Corollary \ref{rang1-bis} to conclude that
  $D'$ is numerically equivalent to a multiple of $H$. Then because $X'$
  is simply connected, this implies the result in the Picard group of
  $X'$.
\end{preu}

\vs 0.2 cm

We finish the proof of Theorem \ref{main}. Because $X'$ is simply
connected, we know that the Picard group is free of finite type. Let
us now prove that its rank is at most
$\frac{1}{2}(h^{n-1}(Y)-\sigma_Y)+1.$ Indeed, let $(D_i)_{0\leq i\leq
  k}$ some divisors on $X'$ independant in $\pic(X')$ with $D_0=H$.
Let us consider the conditions in Corollary \ref{multi}. These
conditions are $x_u\a_u-x_{u^*}\a_{u^*}=0$ for $u\in U$ and in fact
there are only
${\frac{1}{2}h^{n-1}_{hyp}(Y)=\frac{1}{2}(h^{n-1}(Y)-
  \sigma_Y)}$
  such conditions because this condition is trivial if $u=u^*$ and is
  the same for $u$ and $u^*$. Now if $k$ is strictly bigger than the
  preceding number of conditions, then there exists an element $D'$ in
  the linear span of the familly $(D_i)_{1\leq i\leq k}$ satisfying
  these conditions. But then $D'$ as to be a multiple of $H=D_0$, a
  contradiction. 

\section{Applications}
\label{debarre}

Denote by $Y=\G(p,m)$ the Grassmannian variety of $p$-dimensional
vector subspaces in a $m$-dimensional vector space. 
Thanks to \cite[Corollaire 7.4]{De} and \cite[Th{\'e}or{\`e}me 8.1 and
Corollaire 8.3]{De}, we obtain\footnote{The first part
  of this result generalises the celebrated Fulton-Hansen Theorem
  \cite{FH} while the second generalises Bertini's Theorem for
  Grassmannians.}

\vs 0.2 cm

\begin{e-proposition} 
\label{deb}
(\i) Let $X$ be a cumbersome irreducible subvariety in $\G(p,m)$ of
  dimension $n'$ such that $2n'\geq N+r$ (or equivalentely $n'\geq
  n-1+\frac{r}{2}$) then $X$ is simply connected.

\label{irred}
(\i\i) Let $X$ be a cumbersome irreducible variety of dimension $n'\geq
  n-1+\delta$ in $Y$, then the intersection of $X$ with a general
  translate of a Schubert variety of dimension $n$ is irreducible.
\end{e-proposition}

\vs 0.2 cm

%
%
%
These propositions complete the proof of Corollaries \ref{gras} and
\ref{gras-lin}\footnote{The signature $\sigma_Y$ is easy to compute in
  the last case, see \cite{Sl}.}

\vs 0.2 cm

We may not hope that the method proposed by E. Arrondo and
J. Caravantes will lead to a better bound on the rank of the Picard
group because the decomposition into sums of squares shows that some
negative terms appear. More precisely we have the following

\vs 0.2 cm

\begin{e-proposition}
\label{minu}
If $Y$ is cominuscule, for any square of minor appearing $q''$ the
coefficient in $q'$ of this square is smaller than its coefficient in
$q''$.
\end{e-proposition}

\vs 0.2 cm

\begin{preu}
If $t$ is such that $c_{t,H}^uc_{t,H}^{u^*}\neq0$ then the quiver of
$t$ (see \cite{Pe} for more on these quivers) differs from those of
$u$ and $u^*$ by one vertex and has to be the union of these
quivers. In particular $t$ is unique. Furthermore we get $d_t\geq
x_u+x_{u^*}\geq2a_t$ thus $\frac{a_t}{d_t}\leq\frac{1}{2}<1$.
\end{preu}
%
%
%
%
%
%
%
%
%
%
For example, consider the case $Y=Q_{2(n-1)}$ is a smooth quadric. Then
our result together with Fulton-Hansen connectivity Theorem and
Bertini Theorem leads to the:

\vs 0.2 cm

  \begin{e-proposition}
    If $X$ is a smooth subvariety of dimension $n'\geq n$ in a smooth
   quadric $Q_{2(n-1)}$ of dimension $2(n-1)$, then $\pic(X)=\Z$ if
   $n$ is odd and $\pic(X)=\Z$ or $\Z^2$ if $n$ is even.
  \end{e-proposition}

\vs 0.2 cm

This has already been observed by E. Arrondo and J. Caravantes, and in
the last case, they give in \cite{arrondo} an example of smooth
subvariety $X$ of dimension $n$ in $Q_{2(n-1)}$ with Picard group
$\Z^2$.


\end{document}